\documentclass[letterpaper, 14pt]{extarticle}
\usepackage[T1]{fontenc}
\usepackage{times}

\usepackage{doc}

\usepackage{url}

\usepackage{float}
\usepackage{graphicx}
\usepackage{epstopdf}
\usepackage{setspace}
\usepackage{amsmath}
\usepackage{amsfonts}
\usepackage{apacite}
\usepackage{tabularx}
\usepackage{multirow}
\usepackage{scrextend}
\usepackage[font={small}, margin=0.5in]{caption}
\usepackage{hanging}

\usepackage{listings}
\usepackage{color}

\definecolor{mygreen}{rgb}{0,0.5,0}
\definecolor{mygray}{rgb}{0.5,0.5,0.5}
\definecolor{mymauve}{rgb}{0.58,0,0.82}
\lstdefinestyle{customc}{
  belowcaptionskip=1\baselineskip,
  breaklines=true,
  frame=L,
  xleftmargin=\parindent,
  showstringspaces=false,
  basicstyle=\footnotesize\ttfamily,
  keywordstyle=\bfseries\color{mygreen},
  identifierstyle=\color{blue},
  numbers=left,
}

\usepackage{tikz}
\usepackage{pgfplots}

\title{The Base Dependent Behavior of Kaprekar's Routine: A Theoretical and Computational Study Revealing New Regularities}
\author{Daniel Hanover}
\date{}

\makeatletter
\let\thetitle\@title
\let\theauthor\@author
\makeatother

\begin{document}
\onehalfspacing

\maketitle
\thispagestyle{empty}
\pagebreak
\centering{\thetitle}
\abstract{
Consider the following process: Take any four-digit number which has at least two distinct digits. Then, rearrange the digits of the original number in ascending and descending order, take these two numbers, and find the difference between the two. Finally, repeat this routine using the difference as the new four-digit number. In 1949, D. R. Kaprekar became the first to discover that this process, known as the Kaprekar Routine, would always yield 6174 within 7 iterations. Since this number remains unchanged after an application of the Kaprekar Routine, it became known as Kaprekar's Constant. Previous works have shown that the only base 10 Kaprekar's Constants are 495 and 6174, the 3-digit and 4-digit case. However, little attention has been given to other bases or determining which digit cases and which bases have a Kaprekar's Constant. This paper analyzes the behavior of the Kaprekar Routine in the 3-digit case, deriving an expression for all 3-digit Kaprekar Constants. In addition, the author developed a series of C++ programs to analyze the paths integers followed to their respective Kaprekar's Constant. Surprisingly, it was determined from this program that the most commonly required number of iterations required to reach Kaprekar's Constant for 3-digit integers was consistently 3, regardless of base. When loaded as a matrix, the iteration requirement data demonstrates a precise recurring relationship reminiscent of Pascal's Triangle.
}

\thispagestyle{empty}
\pagebreak
\tableofcontents

\thispagestyle{empty}

\clearpage
\pagenumbering{arabic}
\thispagestyle{empty}
%
%
\section{Introduction}

\subsection{Overview}

An Indian schoolteacher and intense number theory devotee, Dattatreya R. Kaprekar has introduced a collection of ideas to recreational mathematics (O'Connor \& Robertson, 2007). One of his ideas involves an operation he created in 1946 known as the Kaprekar Routine (Weisstein, 2004). This procedure involves taking any positive, integral 4-digit number and rearranging its digits in ascending order and in descending order to form two new 4-digit numbers. The output of the Kaprekar Routine is the absolute value of the difference between these two numbers. Applying Kaprekar`s Routine to an integer, n, can be denoted as K(n). He showed that, within 7 iterations, this process would always yield the number 6174, now known as Kaprekar`s Constant, as long as the initially chosen number had at least two distinct digits. The intent of this project is to explore the behavior of Kaprekar's Constant in the 3-digit space. Moreover, it leverages modern computational power to observe unexpected properties of Kaprekar tree systems and produces a novel observation on the Kaprekar iteration distribution graph. 
Although 6174 is technically the only number that has earned the title Kaprekar's Constant, this paper refers to any n-digit number in base b as a Kaprekar Constant if and only if that number is a fixed point under the Kaprekar Routine and all nontrivial n-digit numbers in base b are transformed to that number after a finite number of applications of the Kaprekar Routine. Also, in this paper, a bar over a line of algebra denotes that each expression in that line represents a digit and that the whole line is one number. For example, the expression:
\begin{equation*} \overline{(x+3)(x)(1)}_{10} \end{equation*}
is the equivalent of:
\begin{equation*} 10^2*(x+3)+10^1*x+10^0*1 \end{equation*}

\subsection{Previous Studies}

Kaprekar`s discovery, announced in 1949 at the Madras Mathematical Conference, has piqued the interest of many mathematicians ever since. As a result, the mathematics community has been researching and developing new data on the subject for decades. For instance, Prichett \emph{et al.}, in 1981, determined that the only Kaprekar's Constants in base 10 are 6174 and 495. This was done by showing that for any number of digits greater than 4, there exists a 2-cycle; a cycle occurs when all numbers of a base b with n digits, after several applications of the Kaprekar Routine, cycle through a fixed set of at least 2 n-digit numbers in base b (Jones, 2008). In addition, Ludington concluded in 1979 that for any base, there are a finite number of Kaprekar's Constants (Walden, 2004). 

In 1978, Hasse attempted to find a Kaprekar's Constant for 2 digits. Although the only Kaprekar's Constant found was 01 in base 2, an interesting structure of Kaprekar cycles was found throughout the various bases. This structure is displayed in Figure 1 as a flow chart, wherein every 3 digit number of base 10 is in a box in the nth column (n being the number of times the Kaprekar Routine must be applied before the number reaches a Kaprekar's Constant), and edges connect boxes to other boxes in the next column (n-1th column, since Kaprekar's Constant was just applied once); the second box contains numbers representing the output of the Kaprekar Routine's application. Since a Kaprekar's Constant exists for 3-digit numbers in base 10, all the edges eventually lead to 495; otherwise, some edges would lead to one 3-digit number while others would transpire at another 3-digit number. Each edge is labeled with an integer representing the number of operations flowing along that path (Walden, 2004).

\begin{figure}
\centering
\includegraphics[width=4in]{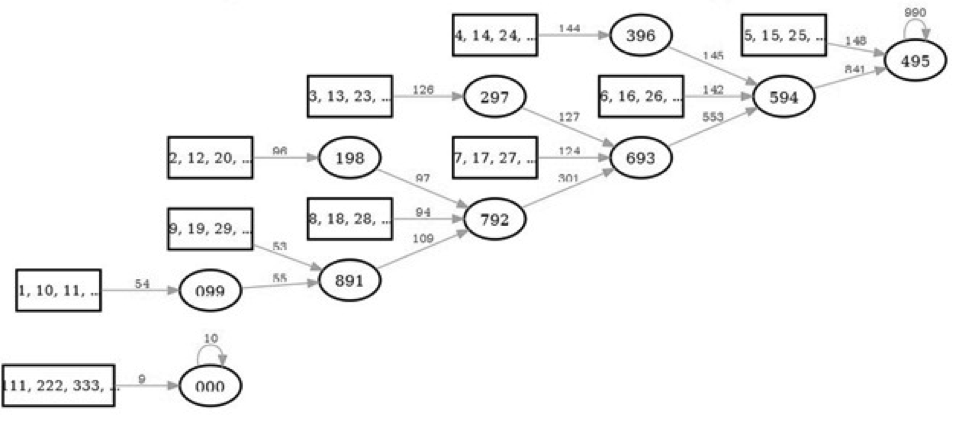} 

\caption{Graph of Kaprekar Routine's Application to 3-digit numbers in base 10. Each edge is labeled with the quantity of transformations flowing through that path. For instance, the square "6, 16, 26, ..." node contains 142 distinct integers, all of which will transform into 594 when the Kaprekar Routine is applied. Adapted from Unknown Creator (15 June 2011). Kaprekar Process for Three Digits. Retrieved from http://en.wikipedia.org/wiki/File:KaprekarRoutineFlowGraph495.svg}
\end{figure}

In 1978, Hasse and Prichett also found that the only bases which had a 4-digit Kaprekar's Constant are $b = 5$ and $b = 4^n*10$. In the former case, the Kaprekar's Constant is 3032. In the latter, the 4-digit Kaprekar's Constant is

\begin{equation*} \overline{(6 * 4^n)(2 * 4^n  -1)(8 * 4^n  -1)(4 * 4^n)} \end{equation*}

$4^n*10$ equals ten when n = 0. Therefore, the 4-digit Kaprekar's Constant of base ten, 6174, would be expected to follow this pattern. Clearly, it does:

\begin{equation*} 6 = 6 * 4^0 \end{equation*}
\begin{equation*} 1 = 2 * 4^0 - 1 \end{equation*}
\begin{equation*} 7 = 8 * 4^0 - 1 \end{equation*}
\begin{equation*} 4 = 4 * 4^0 \end{equation*}

Lastly, Prichett proved that the only bases that have a 5-digit Kaprekar's Constant are those congruent to 3 (mod 6), except for 9 (Walden, 2004). The Kaprekar's Constant for any base $b$, where $r$ is an integer and  $r=(b-3)/6$ can be represented by: 

\begin{equation*} \overline{(4r+2)(2r)(6r+2)(2r+1)} \end{equation*}

While much research has been conducted on Kaprekar's Constant, there are still substantial gaps of information in mathematicians' knowledge. For example, while many researchers have experimented with Kaprekar's Constant in various bases, none have derived a formula which can calculate a Kaprekar's Constant given only its number of digits and its base. Plus, little work has been done on Kaprekar Cycles, which occur when the Kaprekar Routine cycles through a fixed set of at least 2 distinct integers (Jones, 2008). Lastly, some in the mathematical community suspect that there exists deeper meaning to the Kaprekar anomaly; yet, thus far, none have been discovered (Nishiyama, 2006). This project aims to derive a formula that can predict the Kaprekar's Constant for 3-digit numbers in any base. The intent of this paper is also to use computer automation to uncover potential patterns in the Kaprekar Tree for 3-digit numbers in various base systems. 

\section{Methodology}

The basis for many of the discoveries that will be mentioned later in this paper lies in a simplification of the Kaprekar Routine. Consider the Kaprekar Routine expressed in algebraic terms (the number $n$ is expressed in base $b$ with $k+1$ digits):

\begin{equation} n = (\overline{X_{0}X_{1}X_{2}...X_{k}})_{b}, (0 \leq X_k \leq ... \leq X_1 \leq X_0 < b),(X_0,X_1,...,X_k \in \mathbb{Z}) \end{equation}

\begin{equation} n=b^kX_0+b^{k-1}X_1+b^{k-2}X_2 ... + b^0X_k \end{equation}

\begin{equation} K(n)=(b^kX_0+b^{k-1}X_1+b^{k-2}X_2+...) - (b^kX_k+b^{k-1}X_{k-1}+b^{k-2}X_{k-2}+...) \end{equation}

\begin{equation} K(n)=(b^{k-1})(X_0-X_k)+(b^{k-1}-b)(X_1-X_{k-1} )...+(b^{ \lceil (k+1)/2 \rceil} - b^{ \lceil (k+1)/2 \rceil -1} ) (X_{ \lceil (k+1)/2 \rceil -1} - X_{ \lceil (k+1)/2 \rceil }) \end{equation}

Note that the ordering of n's digits does not affect the derivation above. Given any integer $n$ in base $b$ with $k+1$ digits, $X_0$ represents the value of its largest integer, $X_1$ of the next largest, and so on. 
From equation 2.4, it can be concluded that $K(n)$ is always divisible by $(b-1)$, where $b$ is $n$'s base. This is clear because each expression has a factor in the form $b^x-b^y$, where $x>y$. Therefore, $b^y$ can be factored out from each of these factors to get $b^y (b^{x-y}-1)$. Finally, using the Geometric Formula, this simplifies to $b^y (b-1)(1+b+b^2...+b^{x-y-1})$. Since this exists in each expression, $K(n)$ must be divisible by $(b-1)$. Since equation 2.4 makes no assumptions regarding number of digits, base, or ordering of digits, it becomes quite clear that whenever the Kaprekar Routine is applied to an integer, its output will be divisible by the integer's base minus one. While this theorem was discovered and proven before the writing of this paper, and is considered general knowledge, the proof above is novel.

\subsection{Relating to 3-Digit Integers}

An interesting property can be observed with all 3-digit numbers: $K(x_b)$, where x is a 3-digit integer in base b, will always resolve to $\overline{(D-1)(b-1)(b-D)}$, where D is the difference between the largest and smallest digits of x. The author has developed a novel proof that this property holds true for all cases. It is unknown whether it has been proven before, though the property has been observed (Eldridge \& Sagong, 1988). To start, simplify Equation 4 such that it represents only 3-digit integers:

\begin{equation} K(n)=(b^2 - 1)(X_0 - X_2) \end{equation}

In this equation, $(X_0-X_2)$ is equal to D. Using this in Equation 5:

\begin{equation} K(n)=b^2D-D=\overline{D00} - D \end{equation}

Since $D \leq b-1$, the number $\overline{D00}$ minus D will have $D-1$ as the leftmost digit and $b-1$ as the middle digit (because D is not large enough to lower $K(n)$ below $b-1$); the last digit, therefore, must be  $b-D$, since the ones place was 0, before being subtracted by D. In conclusion, when the Kaprekar Routine is applied to any 3-digit number in base b, the result will be $\overline{(D-1)(b-1)(b-D)}$, where D is the difference between the original number's largest and smallest digits.

\subsection{3-Digit Kaprekar's Constants}

In base 10, the 3-digit Kaprekar's Constant is 495. This section of the Methodology aims to determine a formula that will calculate the digits of the Kaprekar's Constant in any base, b.
Finding such a formula for 3-digit integers is relatively simple compared to integers of 5 or even 4 digits. To do so, start with the definition of a Kaprekar's Constant applied to a 3-digit integer found in the previous section:

\begin{equation} K(n)=\overline{(D-1)(b-1)(b-D)} \end{equation}

Define $K^x(n)$ as the Kaprekar function applied to an integer, $n$, $x$ times. Using this notation, it can be said that there are two possibilities for $K^2(n)$:

$D > b/2$

If this is the case, $(D-1)>(b-D)$. This is mainly important because in determining the digits of $K^2 (n)$, the D of $K^1 (n)$ must be known. The difference between the largest and smallest digits of $K(n)$ will be:

\begin{equation} D_{new}=(b-1)-(b-D)=D-1 \end{equation}

Note that $b-1$ is used because, regardless of what base is being used, $b-1$ is the greatest possible value any individual digit can have. Therefore, it must be the largest digit in  $K^2 (n)$. As outlined in the beginning of this paper, for an integer, $n$, to be a Kaprekar's Constant, $K(n)$ must equal $n$. Therefore, $K(n)=K^2 (n)$, if and only if $K(n)$ is Kaprekar's Constant. Using $D_{new}$ as determined above:

\begin{equation} K^2 (n)=\overline{(D_{new}-1)(b-1)(b-D_{new})} = \overline{(D-2)(b-1)(b+D-1)̅} \end{equation}

Comparing $K^2 (n)$ and $K(n)$, it is clear that they can be equal if and only if $D-1=D-2$ and $b-D=b+D-1$. The former is clearly an impossibility; therefore, no Kaprekar's Constant exists when $D>b/2$. 

$D \leq b/2$

If $D≤b/2$, then $(D-1)<(b-D)$. In this case, since $(D-1)<(b-D)$, the difference between the largest and smallest digits of $K^1 (n)$ will be:

\begin{equation} D_{new}=(b-1)-(D-1)=b-D \end{equation}

Following with the same logic presented in the last scenario:

\begin{equation} K^2 (n)= \overline{(D_{new}-1)(b-1)(b-D_{new})} = \overline{(b-D-1)(b-1)(D)} \end{equation}

Comparing $K^2 (n)$ and $K(n)$, it is clear that they can be equal if and only if $D=b-D$ and $(D-1)=(b-D-1)$. Both of these requirements are fulfilled when $b/2=D$. Plugging this in to $K(n)$ and $K^2(n)$:

\begin{equation} K(n)=K^2 (n)=\overline{(b/2-1)(b-1)(b-b/2)} \end{equation}

Therefore, since K(n) was assumed to be a Kaprekar Constant, the number $\overline{(b/2-1)(b-1)(b-b/2)}$ fulfills the first requirement for being a Kaprekar's Constant. 

Nonetheless, to be a true Kaprekar's Constant, it must also be true that every single three digit number in base $b$ (excluding those having all digits identical) will, after repeated iterations of the Kaprekar Routine, transform into this number. Since, in the above two scenarios, only one number was found for which $K(n)=n$, it is not possible for there to be other numbers like it. Therefore, the only way for it not to be a Kaprekar's Constant is if a loop exists, which would mean $K(l)=K^x (l)$ for any integer $x$ and $l$, besides the proposed Kaprekar's Constant. Although the rest of this proof is an adapted version of Eldridge and Sagong's proof, it differs significantly from their version.

Since all numbers are in the form $\overline{(D-1)(b-1)(b-D)}$ after at least one Kaprekar transformation, it can be said that the 1st and 3rd digits always add up to $b-1$. Now, going back to Equation 5, 

\begin{equation} K(n)=(b^2-1)(X_0-X_2) \end{equation}

In this equation, $X_0$ represents the largest digit of n, and $X_2$ represents the smallest one. Since the above holds true for $K(n)$, it also holds true for $K^2 (n)$, except $X_0$ and $X_2$ represent digits of $K(n)$. Since $K(n)$'s largest digit is always $b-1$, $X_0=b-1$. In addition, since the 1st and 3rd digits of any integer which is the result of a Kaprekar Routine's application (such as this integer) add up to $b-1$, and since the second digit is always b-1, it can be concluded that $X_0-X_2=X_1$, where $X_1$ is the digit of $K(n)$ less than $X_0$ and greater than $X_2$. Therefore, $X_1$ could be considered to be D in equation 2.7. This means equation 2.7 can be simplified to: 

\begin{equation} K^2 (n)= \overline{(X_1-1)(b-1)(b-X_1)} \end{equation}

Now, let's compare the digits of $K^2 (n)$ with $K(n)$. The 1st and 3rd digits of $K(n)$ were $X_1$ and $X_2$, because $X_0$, the largest digit, was in the middle. $X_2$ was equal to $b-1-X_1$. Now, $K^2 (n)$ has $X_1-1$ for its 1st digit and $b-X_1$ for its second. Since $X_1-1$ is one less than $X_1$, and since $b-X_1$ is one greater than $b-1-X_1$, it can be concluded that with any Kaprekar Operation on 3-digit integers, one of the outer digits increments and the other decrements. Plus, since $X_1$ is the larger outer digit in $K(n)$ and it is the one being decremented, it can be determined that with every Kaprekar operation, the difference between the outer digits decreases by 2 (This statement, though mentioned in other works, has always been assumed as a Lemma, and never actually proven prior to this paper, as far as the author is aware). For example, take the number 297. The difference between its outer digits is 7 - 2 = 5. By applying the Kaprekar Routine, it is transformed to 972 - 279 = 693. The difference between 693's outer digits is 6 - 3 = 3. 3 is 5 decreased by 2. Thus, the principle holds. 
This is an extremely important point to make: looking back at the proposed Kaprekar's Constant, the difference between the 1st and 3rd digit is 1. Therefore, on any 3-digit number whose outer digits' difference is odd, repeated Kaprekar Operations cannot lead to a looping since the inner digit is always constant $(b-1)$ and the outer digits' difference is always decreasing until it reaches 1. Therefore, it can be said that any 3-digit integer in a base $b$ will, after repeated applications of Kaprekar's Constant be transformed into a number in the form $\overline{(b/2-1)(b-1)(b/2)}$. The only exception to this rule is when $b$ is even, since then and only then will the outer digits' difference be even. Because this number fulfills the second and final requirement for it to be a Kaprekar's Constant, it can be concluded that it is the one.

\subsection{Computations}

The first program developed for analysis was a C++ code that iterated through all 3- and 4-digit integers, repeatedly applying Kaprekar's Routine until the respective Kaprekar's Constant was found, and then moving on to the next integer. The program was adapted to also track the number of iterations of the Kaprekar Routine required for each integer to reach its Kaprekar's Constant. However, the program, when initially conceived, was wildly inefficient: the runtime to crunch through 4-digit integers was over an hour, ten times that of crunching the 3-digit integers. Predictably, processing larger orders of magnitudes became impractical without a refinement of the algorithm. 
The clear issue was that the algorithm was testing every possible integer case, treating every integer separately. However, mathematically, it is possible to group many of these integers together. Walden (2004) devised the following encoding schema:

Rewritten, the Kaprekar Routine applied to an N-digit number where the digits $d_{N-1} \geq d_{n-2} \geq ... \geq d_0$ is simply:

\begin{equation*}
\begin{matrix}
    \space & d_{N - 1}    & d_{N - 2} & \dots & d_{1} & d_{0} \\
    - & d_{0}       & d_{1} & \dots & d_{N-2} & d_{N-1} \\ \hline
\end{matrix}
\end{equation*}

However, this can also be written as: 

\begin{equation*}
\begin{matrix}
    \space & d_{N - 1}-d_0    & d_{N - 2}-d_1 & \dots & 0 & 0 \\
    - & 0       & 0 & \dots & d_1 - d_{N-2} & d_0 - d_{N-1} \\ \hline
\end{matrix}
\end{equation*}

Clearly, it is redundant to include the latter 0's, and any number can be represented with the differentials of its respective digits. Thus, when referring to Kaprekar's Constant, any N-digit number can be referred to using only $\lfloor N/2 \rfloor$ digits: For instance, 8991 can be represented as $<9-1|9-8> \rightarrow <8|1>$. The Kaprekar Routine applied to any number that fits this notation would have the same output as if it were applied to 8991.

Walden (2004) also developed a caching scheme based on this mathematical property. However, this caching scheme was insufficient for the needs of this study because it had no mechanism for storing how far each node was from its Kaprekar Constant. Furthermore, the indexing system of Walden's caching system did not relate each node to its position in the underlying Kaprekar Tree, making it difficult to track the iteration requirements of each node.
With this in mind, the author developed a caching scheme, distinct from Walden's but still taking advantage of his encoding schema, to incorporate into the program, reducing runtime by a factor of 100. Specifically, a recursive node-linked C structure was built (see Appendix A). The program iterated through all encodings from $<0,0,0…>$ to $<b-1,b-1,b-1>$, filling the recursive structure as new node pathways were discovered. Each node also stored the number of iterations it would need to reach the Kaprekar's Constant node. This meant that rather than reaching the Kaprekar's Constant node for every integer case, the program merely had to reach a node it had already solved, and it would be able to use that node's cached iteration value. This program was used to map out the 3-digit integer case for all even bases from 2 to 56.

\section{Results \& Discussion}

In summary, a basic background of Kaprekar's Constant and its properties were discussed as well as a variety of proofs. It was proven that, whenever the Kaprekar Routine is applied to an integer, the result would be divisible by that number's base minus 1. In addition, the result when the Kaprekar Routine is applied to a 3-digit integer in any base was found in terms of that number's base as well as the difference between its largest and smallest digits. Thirdly, it can now be said that whenever the Kaprekar Routine is applied to a 3-digit integer, if that integer was the result of a previous Kaprekar Routine, then the difference between the outer digits of the 3-digit integer will be decreased by two. Finally, the Kaprekar's Constant for all 3-digit integers in an even base was found in terms of that number's base.
Some conclusions were also derived from the C++ program created to iterate through all 4-digit and 3-digit integers and test each of them with the Kaprekar Routine. The program was modified to create a list at the end of its execution mapping numbers to the quantity of integers that required that many iterations of the Kaprekar Routine to reach a Kaprekar Constant. The results are displayed in Figures 2 and 3.

\begin{figure}[H]
\centering
\includegraphics[width=4in]{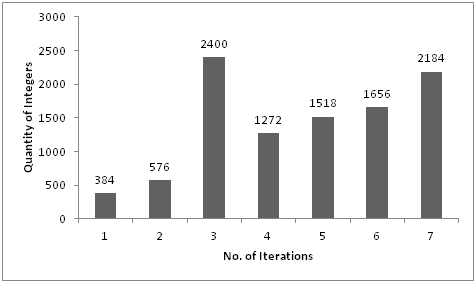} 

\caption{Kaprekar Routine Iterations: 4-Digit. Demonstrates the relationship between Kaprekar Routine iterations and the number of base 10,  4-digit integers which require them to reach a Kaprekar Constant. Graphic by author}
\end{figure}

\begin{figure}[H]
\centering
\includegraphics[width=4in]{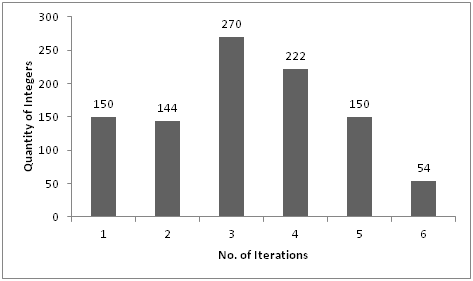} 

\caption{Kaprekar Routine Iterations: 3-Digit. Demonstrates the relationship between Kaprekar Routine iterations and the number of base 10, 3-digit integers which require them to reach a Kaprekar Constant. Graphic by author}
\end{figure}

Note that, in both cases, most integers were 3 iterations of the Kaprekar Routine away from the Kaprekar Constant. Taking advantage of the caching scheme developed earlier, iteration requirement graphs were generated for all 3-digit integers from bases 2 to 56, and a subset of the results are shown in Table 1. Figure 4 visualizes this data and depicts a fascinating pattern in the data. Strangely, regardless of base, the most commonly required number of iterations is consistently 3. It is also curious to observe that, despite different bases having completely different Kaprekar tree structures, they maintain nearly parallel iteration distribution curves. 

\begin{table}[H]
\centering
\caption{3-Digit Kaprekar Iteration Distribution: Bases 2-32, Tabular Form. A number of patterns between cells are readily apparent, much like Pascal's Triangle}
\small
\begin{tabularx}{6.5in}{cc|XXXXXXXX}
\multicolumn{9}{c}{\small Number of Iterations} \\
\multirow{15}{*}{\rotatebox[origin=c]{90}{\small Base}} &    & 1    & 2    & 3    & 4    & 5    & 6    & 7    & 8    \\ \cline{2-10}
 & 2  & 6    & 0    & 0    & 0    & 0    & 0    & 0    & 0    \\
 & 4  & 24   & 18   & 18   & 0    & 0    & 0    & 0    & 0    \\
 & 6  & 54   & 48   & 78   & 30   & 0    & 0    & 0    & 0    \\
 & 8  & 96   & 90   & 162  & 114  & 42   & 0    & 0    & 0    \\
 & 10 & 150  & 144  & 270  & 222  & 150  & 54   & 0    & 0    \\
 & 12 & 216  & 210  & 402  & 354  & 282  & 186  & 66   & 0    \\
 & 14 & 294  & 288  & 558  & 510  & 438  & 342  & 222  & 78   \\
 & 16 & 384  & 378  & 738  & 690  & 618  & 522  & 402  & 258  \\
 & 18 & 486  & 480  & 942  & 894  & 822  & 726  & 606  & 462  \\
 & 20 & 600  & 594  & 1170 & 1122 & 1050 & 954  & 834  & 690  \\
 & 22 & 726  & 720  & 1422 & 1374 & 1302 & 1206 & 1086 & 942  \\
 & 24 & 864  & 858  & 1698 & 1650 & 1578 & 1482 & 1362 & 1218 \\
 & 26 & 1014 & 1008 & 1998 & 1950 & 1878 & 1782 & 1662 & 1518 \\
 & 28 & 1176 & 1170 & 2322 & 2274 & 2202 & 2106 & 1986 & 1842

\end{tabularx}

\end{table}

\begin{table}[H]
\centering
\caption{Relationships between numbers in Table 1. Columns rise quadratically, whereas diagonals follow a linear path. Note also that there is a significant shift in behavior between columns 2 and 3.}
\includegraphics[width=4in]{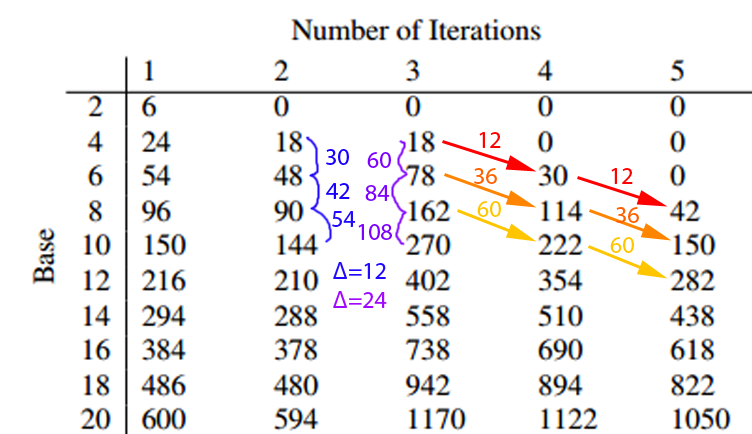} 
\end{table}

Table 1 is the numeric representation of the Kaprekar Iteration Requirement tree. Each cell represents the quantity of integers in base R that require C iterations to converge on Kaprekar's Constant. Every integer in Table 1 is divisible by 6. This result is not particularly surprising, however, since every 3-digit integer with unique digits has 5 permutations, all of which have exactly the same iteration requirements. 

A deeper analysis reveals less trivial regularities. Notice that the maximum number of iterations required consistently increases by 1 each time the base increments. Moreover, the number of integers that require the largest iteration count consistently increases by 12 with each new row in the table. For instance, in base 4, 18 numbers require 3 iterations to converge, whereas in base 6 for 4 iterations, this number rises by 12 to 30. Each other diagonal in the table follows a constant slope as well, with the slope increasing by 24 in each diagonal further from the edge of the triangle. 

There also seems to be a fundamental change as the table moves past the second data column. For instance, the diagonal slopes do not hold true in the first two columns. Moreover, in columns 1 and 2, the difference between each cell and the cell below it increases by 12 with each row. However, from column 3 onward, this delta doubles to 24, as mentioned in Table 2. Considering these irregularities, as well as the jump in Figure 4 that occurs at each line's 3rd data point, it becomes clear that some peculiar effect is transforming the behavior of these Kaprekar Systems as integers 3 layers away from Kaprekar's Constant are considered. Note, also, that as a direct result of the patterns observed above, each column of the table follows a perfect quadratic curve ($y=6x^2+b, y=12x^2+b$), and each diagonal follows a perfect linear path (with the outer diagonal following $y=12x$). Nonetheless, most of these observations have no explanation, and there is more work necessary to link this peculiar graph with the the current understanding of the underlying Kaprekar trees.

\begin{figure}
\centering
\includegraphics[width=4in]{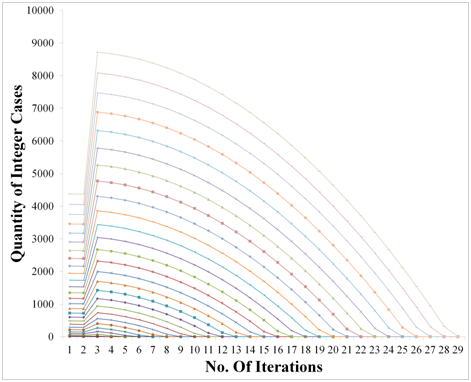} 

\caption{3-Digit Kaprekar Iteration Distribution: Bases 2-56. Quantity of Integers as a function of the number of Kaprekar Iterations those integers need to reach Kaprekar's Constant. This graph presents data from 3-digit integers, each line representing a different even base between 2 and 56. Note that, in all base systems above 4, the most commonly required number of iterations is 3. Graphic by author}
\end{figure}

Lastly, Walden, 2004 developed a program to determine the Kaprekar tree of a given base with a given number of digits. This program was modified by the author to observe the behavior of Kaprekar Cycles. Cycles occur when $K(l)=K^x (l)$ and $K(l) \neq K^{x-i} (l)$, for any integers $x$, $l$, and $i$, where $i$ is less than $x$. In other words, a Kaprekar Cycle is a loop in the Kaprekar tree. Specifically, the program generated a line graph, shown in Figure 5, plotting the quantity of n-cycles as a function of n, where an n-cycle is a Kaprekar Cycle with $x=n$, for all positive integers between base 2 and 36 and comprising of between 2 and 9 digits, inclusive. Notice that while the most common Kaprekar Cycle length is 1, the second-most common cycle length is, surprisingly, 5. There is no known mathematical explanation of why 5-cycles are particularly common.

\begin{figure}
\centering
\includegraphics[width=4in]{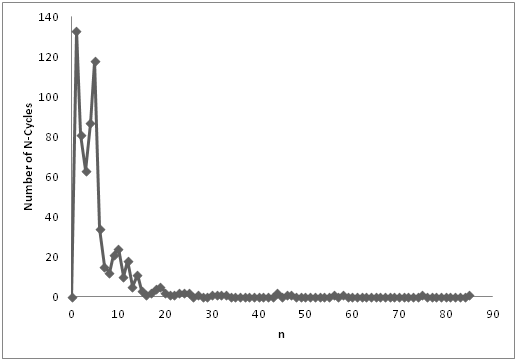} 

\caption{Distribution of Cycles. Demonstrates the relationship between N and the number of N-cycles that appear in the Kaprekar tree.}
\end{figure}

\section{Continuing Research}

Although a lot of progress has been made in this paper, there is still much mystery surrounding the number theory behind Kaprekar's Constant. For example, while 3-digit integers were the prime focus of this paper, research is being done on how 5-digit and 9-digit numbers relate to Kaprekar's Constant. Also, the unexpected result of graphing the iteration distribution for the various bases with 3-digit integers was not fully explored in this paper, and a more detailed analysis will be necessary. Lastly, by employing customized computer programs for data testing and analysis, new Lemmas are currently being developed that will aid in the creation of deeper, more complex laws describing the behavior of Kaprekar's Constant.

\section{References}
\begin{sloppypar}
  Eldridge, K. E., \& Sagong, S. (1988).
  The Determination of Kaprekar Convergence and Loop Convergence of All Three-Digit Numbers.
  \emph{In The Determination of Kaprekar Convergence and Loop Convergence of All Three-Digit Numbers} (pp. 105-112). 
  Amer. Math. Monthly no. 2.
  
  Jones, E. (2008, July). \emph{The Kaprekar Routine}. Retrieved November 5, 2011, from University of Nebraska-Lincoln: http://scimath.unl.edu/MIM/files/MATExamFiles/Jones\_MATpaper\_Final\_071408\_DMc.pdf
  
  Nishiyama, Y. (2006, March 1). \emph{Mysterious number 6174}. Retrieved November 5, 2011, from Plus Magazine: http://plus.maths.org/content/os/issue38/features/nishiyama/index

  O'Connor, J. J., \& Robertson, E. F. (2007, August). \emph{Dattatreya Ramachandra Kaprekar}. Retrieved November 5, 2011, from The MacTutor History of Mathematics archive: http://www-history.mcs.st-andrews.ac.uk/Biographies/Kaprekar.html
  
  Walden, B. L. (2004, April 7). \emph{Searching for Kaprekar's Constants: Algorithms and Results}. Retrieved December 10, 2011, from Hindawi Publishing Corporation: http://downloads.hindawi.com/journals/ijmms/2005/802780.pdf

  Weisstein, E. W. (2004, November 11). \emph{Kaprekar Routine}. Retrieved November 5, 2011, from MathWorld--A Wolfram Web Resource: http://mathworld.wolfram.com/KaprekarRoutine.html 
  \end{sloppypar}
\section{Appendix A}
To prove that the Kaprekar's Constant of 4-digit numbers in base 10 is 6174, a C++ program was written to iterate through all numbers between 1000 and 9999 (excluding multiples of 1111) and apply the Kaprekar Routine on each integer repeatedly, until 6174 was reached. The result of running the program confirmed D. R. Kaprekar's observation; every integer reached 6174 within 7 applications of the Kaprekar Routine. The program was then modified, incorporating an encoding scheme similar to one employed by Walden, 2004 to vastly improve running time. It was also modified to record the number of iterations required for each input to converge on Kaprekar's Constant. The important excerpts of the source code of the program are show below:
\begin{singlespace} \noindent
\begin{lstlisting}[language=C++, style=customc]
#include <iostream>
#include <iomanip>
#include<math.h>
#include<stdlib.h>

using namespace std;

int searchaway(node *p, ArrayPtr key, const int& N) {
	int n=N/2, c=compare(p->ptr,key, n), j, k;
	if (c==-1)  
	{
		if (p->l==NULL)
		{
			cout << "error" << endl;
		}
		else {
			searchaway(p->l,key, N);
		}
	}
	else if (c==1) 
	{
		if (p->r==NULL)
		{
			cout << "error" << endl;
		}
		else {
			searchaway(p->r,key, N);
		}
	}
	else if (c==0)
	{
		return p->away;   
	} 
}
node *searchadd(node *p, ArrayPtr key, bool& found, bool& newl, 
				const int& N, int& numpass, int& numpos, const int& top)
{ 
	//adds node to p; omitted for brevity 
}  

int main(int argc, char *argv[])
{
	//initialization
	int  i, n, N=atoi(argv[1]),  
		base, topdig, numpass=1, numpos=0;
	ArrayPtr Q, R;
	node *root;
bool found=false, newl=false;

	if (argc == 5)
	{
		n=N/2;
		for(base=atoi(argv[2]); base<=atoi(argv[3]); base+=atoi(argv[4]))
		{
			cout << base << endl << endl;
			topdig=base-1;
			Q = new int[n];
			if (Q==NULL)
			{
				cout << "Insufficient storage.\n";
				exit(1);
			}
			R = new int[n]; 
			if (R==NULL)
			{
				cout << "Insufficient storage.\n";
				exit(1);
			} 
			Q[0]=0;
			for(i=0;i<n;i++)
				Q[i]=0;  

			root = new node;
			if (root==NULL)
			{
				cout << "Insufficient storage.\n";
				exit(1);
			}
			root->l=NULL;
			root->r=NULL;
			root->ptr=Q;
			root->pass = numpass;
			root->pos = numpos;
			root->away = -1;

			node** res = new node*[100];
			int c;
			//for(c=0; c<100; ++c)
			//	res[c] = new int[n];
			c=0;
			Q=nxt(Q,topdig,N);   
			while(Q[0]!=0)
			{
				c++;
				found=false; newl=false;
				res[c] = searchadd(root, Q, found, newl, N, numpass, numpos, topdig);
				if (!(found)) {
					R=Q;
				}
				while(!(found))
				{
					c++;
					R=succ(R,topdig,N);
					res[c] = searchadd(root, R, found, newl, N, numpass, numpos, topdig);
				}

				int aw = 0;
				for (int k = c; k >= 1; k--)
				{
					if (res[k]->away != -1)
					{
						aw = res[k]->away;
						continue;
					}
					aw++;
					res[k]->away = aw;
}
				c=0;
				Q=nxt(Q,topdig,N);  
			}

			delete [] Q;
			delete [] R;

			int MAX_CYCLES = 150;
			int* a = new int[MAX_CYCLES]; //100 is arbitrary, should be the maximum possible number of cycles
			for(int i = 0; i < MAX_CYCLES; i++) 
				a[i] = 0;

			for(int i = 0; i < pow(base, N); i++)
			{
				int counter = 0;
				int startNum = i;
				ArrayPtr arr = new int[N];
				while(counter < N) {
					arr[counter++] = startNum%base;
					startNum /= base;
				}
				BubbleSortArray(arr, N);
				ArrayPtr Q = new int[n];
				for (int i = 0; i < n; i++)
				{
					Q[i] = arr[N - i - 1] - arr[i];
				}
				int cycles = searchaway(root, Q, N);
				if (cycles >= 0)
				{
					a[cycles]++;
				}
				delete [] arr;
				delete [] Q;
			}

			cout << "Results for base = " << base << ", N = " << N << ": " << endl;
			int highestIndex = 0;
			int highestAmt = 0;
			for (int q = 0; q < MAX_CYCLES; q++)
			{
				if (a[q] > highestAmt)
				{
					highestAmt = a[q];
					highestIndex = q;
				}
				if (a[q] != 0)
				{
					cout << (q+1) << ": " << a[q] << endl;
				}
			}
			trash(root); 
			delete root;
			delete [] res;
			delete [] a;
		}
		return 0;
	} 
}
\end{lstlisting}
\end{singlespace}

\end{document}